\documentclass[11 pt]{amsart}

\usepackage[pagewise]{lineno}
\usepackage[top=1.4in, bottom=1.25in, left=1.25in, right=1.25in]{geometry}

\usepackage{mathrsfs}

\usepackage{hyperref, amsmath, amsthm, amsfonts, amscd, flafter,epsf,amssymb}

\newtheorem{theorem}{Theorem}
\newtheorem{lemma}{Lemma}
 \newtheorem*{thm 1}{Theorem 1}
 
\newcommand{\be}{\begin{equation}}
\newcommand{\ee}{\end{equation}}

\renewcommand{\b}{\beta}

\newcommand{\g}{\gamma}
\newcommand{\s}{\sigma}
\newcommand{\di}{i} 

\numberwithin{equation}{section}
\usepackage{parskip}
\begin{document}

\title{Extreme values of the Riemann zeta function at its critical points in the critical strip}
\author{Shashank Chorge}
\address{Department of Mathematics, University of Rochester, Rochester NY 14627 USA}
 
\maketitle

\begin{abstract}  
We estimate   large  and small values of $|\zeta(\rho')|$, where $\rho'$ runs over  critical points of 
the zeta function in the right half of the critical strip, that is,  the points where $\zeta'(\rho')=0$ and $1/2<\Re \rho'<1$.\\
\end{abstract}

\section{Introduction}

 Assuming the Riemann hypothesis (RH), Littlewood\cite{paper3}  proved that
 \be\label {Litt 1}
 \limsup\limits_{t \rightarrow\infty}\frac{|\zeta(1+\di t)|}{\log\log t}\leq 2e^{C_0},
 \ee
where $C_0$ denotes Euler's constant.
He also proved, unconditionally,  that 
\be\label {Litt 1}
\limsup\limits_{t \rightarrow\infty}\frac{|\zeta(1+\di t)|}{\log\log t}\geq e^{C_0}.
\ee

Gonek and Montgomery
~\cite{paper} obtained similar results sampling zeta at the critical points of the Riemann zeta function to the right of $\Re s=1.$
Let $\rho'=\beta'+\di\gamma'$ denote a typical critical point of the zeta function, that is, a point where $\zeta'(\rho')=0$.
Assuming RH, they showed that for $\beta' \geq 1$,
\be\label {MG 1}
\limsup\limits_{\gamma' \rightarrow\infty}\frac{|\zeta(\rho')|}{\log\log \gamma'}\leq \frac{1}{2}e^{C_0}
\ee
and, unconditionally, for $\beta'>1$, that
\be\label {MG 2}
\limsup\limits_{\gamma' \rightarrow\infty}\frac{|\zeta(\rho')|}{\log\log \gamma'}\geq \frac{1}{4}e^{C_0}.
\ee
Gonek and Montgomery mention that one of their motivations was to answer a question posed by J. G. Thompson  
as to whether, for any large constant $c$, there must always be infinitely
many compact connected components of the level set $|\zeta(s)|=c$.  We see from
\eqref{MG 1} that the answer is yes, for if we take $c$ to be just slightly
less than $|\zeta(\rho')|$,  there will be at least one compact connected component
that passes close to $\rho'$.  Furthermore, Montgomery and Thompson \cite{Mont-Thomp} have shown
that the imaginary parts of points on a connected compact component of a level set
$|\zeta(s)|=c$ will all lie in an interval of the form $[T-C\log T,T+C\log T]$ where
$C$ is an absolute constant.  Thus widely-spaced critical points $\rho'$ will give rise to
 disjoint compact connected components.

In this paper we obtain estimates corresponding to \eqref{MG 1} and  \eqref{MG 2} in the right half of the critical strip.
\begin{theorem}\label{thm 1}
 Assume RH. Let $\sigma_1$ and $\sigma_2$ be fixed with $1/2<\sigma_1<\sigma_2<1.$ 
 If $\rho'= \beta'+\di \gamma'$ is any critical point of the Riemann zeta function with
 $\sigma_1<\beta'<\sigma_2,$ then there is a positive constant $A$ depending on $\s_1$ and $\s_2$ such that
\begin{equation}\label{est1}
\log |\zeta(\rho')|\le A  \ {\frac{(\log \gamma')^{2-2\b'}}{\log \log \gamma'}}.
\end{equation}
\end{theorem}

The estimate in (\ref{est1}) is in fact true with $\rho'$ replaced by any point $s=\sigma+\di t$ with
$\sigma_1<\sigma <\sigma_2$ (see, for example, equation (14.5.2) of Titchmarsh~\cite{Titch}), so it follows from it. 
We nevertheless give a separate proof  since the method we use informs the proof of our second theorem.

\begin{theorem}\label{thm 2}
Let $\sigma_1$ be such that $1/2<\sigma_1<1$, 
Let $\rho'=\beta'+\di \gamma'$ denote a critical point of the Riemann zeta function such that $\sigma_1<\beta'<1$ . Then for any $\epsilon$
and $\epsilon'>0$ and for infinitely many $\rho'$ with $\gamma'\to\infty$, we have unconditionally that
\begin{equation}\label{eqn thm 2}
 \log|\zeta(\rho')|\geq 
 (B(\sigma_1)  -\epsilon' )\frac{(\log\gamma')^{1-\beta'}}{(\log\log\gamma')^{5-3\beta'+\epsilon}},
  \end{equation}
  where
  \be
  B(\sigma_1)  =\frac{(\sigma_1-1/2)^{1-\s_1}\log 2}{(1-\s_1)^2 4^{1-\s_1}}.
  \ee
\end{theorem}

One can ask for similar results 
concerning the small values of $\zeta(1+\di t) $. Assuming RH, Littlewood~\cite{paper4} showed that $$\liminf\limits_{t \rightarrow\infty}|\zeta(1+\di t)|\log\log t\geq \frac{\pi^2}{12}e^{-C_0},$$ and Tichmarsh~\cite{paper5} showed unconditionally that $$\liminf\limits_{t \rightarrow\infty}|\zeta(1+\di t)|\log\log t\leq \frac{\pi^2}{6}e^{-C_0}.$$
Gonek and Montgomery~\cite{paper} obtained corresponding  results at the critical points of the Riemann zeta function to the right of $\Re s=1$.
They showed unconditionally, for $\beta'>1$, that
$$\liminf\limits_{\gamma' \rightarrow\infty}{|\zeta(\rho')|}{\log\log \gamma'}\geq \frac{\pi^2}{3}e^{-C_0},$$ 
and, assuming RH,  that
$$\liminf\limits_{\gamma' \rightarrow\infty}{|\zeta(\rho')|}{\log\log \gamma'}\leq \frac{2\pi^2}{3}e^{-C_0}.$$
We prove the  analogous results at the critical points  $\rho'$  in the right half of the critical strip.
\begin{theorem}
Assume RH. 
  Let $\sigma_1$ and $\sigma_2$ be such that $1/2<\sigma_1<\sigma_2<1.$ 
 If $\rho'= \beta'+\di \gamma'$ is any critical point of the Riemann zeta function with
 $\sigma_1<\beta'<\sigma_2,$
 then for some constant $C$ depending on $\sigma_1$ and $\sigma_2$, we have
 \begin{equation}\label{est3}
   \log |\zeta(\rho')|\geq\frac{-C(\log \gamma')^{2-2\beta'}}{\log \log \gamma'}.  
 \end{equation}
\end{theorem}

Here, as in Theorem~\ref{thm 1}, the estimate in (\ref{est3}) is true even if we replace $\rho'$ with a generic point $s=\sigma+\di t.$
\begin{theorem}
Let $\sigma_1$ be such that $1/2<\sigma_1<1$ and, let $d_1=\sigma_1-1/2.$
Let $\rho'=\beta'+\di \gamma'$ denote a critical point of the Riemann Zeta function such that $\sigma_1<\beta'<1$. Then for any $\epsilon$ and $\epsilon'>0$ and for infinitely many $\rho'$ with $\gamma'\to\infty$, we have unconditionally that
\begin{equation}\label{eq1.9}
 \log|\zeta(\rho')|\leq (-B(\s_1)+\epsilon')\frac{(\log\gamma')^{1-\beta'}}{(\log\log\gamma')^{5-3\beta'+\epsilon}},
\end{equation}
where
\be
  B(\sigma_1)  =\frac{(\sigma_1-1/2)^{1-\s_1}\log 2}{(1-\s_1)^2 4^{1-\s_1}}
  \ee
\end{theorem}
 
\section{Lemmas for the proof of Theorem 1} 
\begin{lemma}\label{lemma1}
	Assume RH. Let $\sigma_1$ and $\sigma_2$  be fixed with $1/2<\sigma_1<\sigma_2<1$.
Then for any $\sigma$  with $\sigma_1\leq \sigma\leq \sigma_2$  and any $x\geq 2$ 
 $$\sum_{ n\leq x} \frac{\Lambda(n)}{n^\sigma}= \frac{x^{1-\sigma}}{1-\sigma} + O(1),$$
 where the implied constant in the $O$-term depends on $\sigma_1$ and $\sigma_2$.
 \end{lemma}
\begin{proof}
Let
$$\psi(x) = \sum_{n\leq x}\Lambda(n),$$
where $\Lambda(n)$ is von Mangoldt's function, namely, $\Lambda(n)=\log p$ if $n$ is a power 
of the prime $p$ and $\Lambda(n)=0$ otherwise.
Assuming the Riemann hypothesis, we have
\be
\psi(x)=x+E(x),
\ee
where $E(x) \ll x^{1/2}\log^2 x$. Thus, using 
Stieltje's integration and integration by parts, we see that
\be
\begin{split}
\sum_{n\leq x} \frac{\Lambda(n)}{n^\sigma}=\int_{1}^{x} \frac{d\psi(y)}{y^\sigma}
=&\int_{1}^{x}  y^{-\sigma} dy +   \int_{1}^{x} \frac{dE(y)}{y^\sigma}\\
=&\frac{x^{1-\sigma}-1}{1-\sigma} + \frac{E(x)}{x^{\sigma}} +\sigma\int_1^x \frac{E(y)}{y^{\sigma+1}}
\end{split}
\ee	
The last two terms are 
$
\ll x^{1/2-\sigma_1}\log^2x \;\big(1+ \tfrac{1}{\sigma_1-1/2} \big)\ll_{\sigma_1} 1. 
$	
Furthermore,   $1/(\sigma-1) \ll_{\sigma_2} 1$. Hence
 $$\sum_{ n\leq x} \frac{\Lambda(n)}{n^\sigma}= \frac{x^{1-\sigma}}{1-\sigma} + O(1)$$
with the implied $O$-term constant depending  on $\sigma_1$ and $\sigma_2$, as asserted.
\end{proof}

In our proofs from now on, we will not keep  track as explicitly as above of the dependence of our $O$-term constants on the 
parameters $\sigma_1$ and $\sigma_2$.

\begin{lemma}\label{lemma2} 
With  the same hypotheses as in Lemma 1, we have
$$\sum_{2\leq n\leq x} \frac{\Lambda(n)}{n^\sigma \log n}=\frac{x^{1-\sigma}}{(1-\sigma)\log x}+O\bigg({\frac{x^{1-\sigma}}{{\log^2 x}}}\bigg).$$\\
Here the implicit $O$-term constant depends on $\sigma_1$ and $\sigma_2.$
\end{lemma}

\begin{proof}
 Define $$S(y) = \sum_{  n\leq y}\frac{\Lambda(n)}{n^\sigma}.$$ From Lemma 1 we find that
\begin{align*}
	  \sum_{ n \leq x} \frac{\Lambda(n)}{n^\sigma \log n} &= \int_{2}^{x} \frac{dS(y)}{\log y}= \frac{S(x)}{\log x}+\int_{2}^{x}\frac{S(y)}{y \log^2 y }dy\\
	  &= \frac{x^{1-\sigma}}{(1-\sigma)\log x} + O\bigg(\frac{1}{\log x}\bigg) +\int_{2}^{x}\frac{y^{1-\sigma}}{y\log^2 y}dy+\int_{2}^{x}O\bigg(\frac{1}{y\log^2 y}\bigg)dy\\
&=\frac{x^{1-\sigma}}{(1-\sigma)\log x}+O\bigg({\frac{x^{1-\sigma}}{{\log^2 x}}}\bigg),
\end{align*}
where the implied constant depends on $\s_1$ and $\s_2$.	
\end{proof}

\begin{lemma}\label{lemma3}
  Let $\sigma_1$ and $\sigma_2$ be such that $1/2<\sigma_1<\sigma_2<1$. Then for $\sigma_1<\sigma<\sigma_2$,
 $$
 \Re\sum_ { n\leq x} \frac{\Lambda(n)}{n^s \log n}\leq\frac{x^{1-\sigma}}{(1-\sigma)\log x} +O\bigg({\frac{x^{1-\sigma}}{\log^2 x}}\bigg),
 $$ 
the implied constant depending on $\s_1$ and $\s_2$. 
\end{lemma}

\begin{proof}
We have
\begin{align*}
\Re\sum_{2\leq n\leq x} \frac{\Lambda(n)}{n^s \log n}\leq \sum_{2\leq n\leq x} \frac{\Lambda(n)}{n^\sigma \log n}=\frac{x^{1-\sigma}}{(1-\sigma)\log x}+O\bigg({\frac{x^{1-\sigma}}{{\log^2 x}}}\bigg)
\end{align*}
by Lemma 2.
\end{proof}

\begin{lemma}\label{log deriv}
Assume RH. Let $\sigma_1$  be such that  $1/2<\sigma_1<1$ and  let $4\leq T\leq t\leq 2T.$
Then for $\sigma_1\leq\sigma<1$,
    $$-\frac{\zeta'}{\zeta}(s)=\sum_{ n\leq \log^2 T} \frac{\Lambda(n)}{n^s}+O(\log^{2(1-\sigma)} T). 
 $$
The implied constant depends at most on $\s_1$.  
\end{lemma}

\begin{proof}
	From a theorem of Montgomery and Vaughan(\cite{classical}, (13.35)), when $s$ is not a root or the pole of the Riemann zeta function and $x, y\geq 2$, we have
	
\begin{align}\label{MV expl frml}
    -\frac{\zeta'}{\zeta}(s)=&-\sum_{n\leq xy} w(n)\frac{\Lambda(n)}{n^s}+\frac{(xy)^{1-s}-x^{1-s}}{(1-s)^2\log y}-\sum_{\rho}\frac{(xy)^{\rho-s}-x^{\rho-s}}{(\rho-s)^2\log y}\\&-\sum_{k=1}^{\infty}\frac{(xy)^{-2k-s}-x^{-2k-s}}{(2k+s)^2\log y},
\end{align}	

where
$$
 w(u)=
\begin{cases}
	1  &\qquad \hbox{if}\qquad 1\leq u\leq x,\\
	1  -\frac{\log(u/x)}{\log y} &\qquad  \hbox{if}\qquad  x<u\leq xy,\\
	 0 &\qquad  \hbox{if}\qquad  u>xy.				
\end{cases}
$$	        	
Substituting $y=2$ and $x=\log^2 T$ in the second term on the right-hand side and taking absolute values, we get
$$\frac{(xy)^{1-s}-x^{1-s}}{(1-s)^2\log y}=O\bigg(\frac{x^{1-\sigma}}{T^2}\bigg)=O\bigg(\frac{\log^{2-2\sigma}T}{T^2}\bigg).$$
The last term on the right-hand side is absolutely bounded, thus  $O(1).$
To simplify the third term on the right-hand side
we assume the Riemann hypothesis and that $\sigma>\sigma_1>1/2.$ 
Substituting $x=\log^2T$ and $y=2$, and then taking absolute values, we find that $$\sum_{\rho}\frac{(xy)^{\rho-s}-x^{\rho-s}}{(\rho-s)^2\log y}=O\bigg(\sum_{\rho}\frac{x^{1/2-\sigma}}{|\rho-\sigma|^2}\bigg)
=O (x^{1/2-\sigma}{\log T} )=O(\log ^{2-2\sigma}T).$$
Here the final two $O$-term constants depend on $\sigma_1.$
Now consider the first term on the right-hand side, which is,
\begin{align*}
	\sum_{n\leq xy} w(n)\frac{\Lambda(n)}{n^s}&=\sum_{n\leq x}\frac{\Lambda(n)}{n^s}+ \sum_{x < n\leq 2x}w(n)\frac{\Lambda(n)}{n^s}=\sum_{n\leq x}\frac{\Lambda(n)}{n^s}+ O\Big(x^{-\sigma}\sum_{x< n\leq 2x}\Lambda(n)\Big)\\
	&=\sum_{n\leq x}\frac{\Lambda(n)}{n^s}+ O(x^{1-\sigma})=\sum_{n\leq \log^2 T}\frac{\Lambda(n)}{n^s}+ O(\log^{2-2\sigma}T),
	\end{align*}
	where here the implied constant is absolute. 
Combining all our results, we obtain
$$-\frac{\zeta'}{\zeta}(s)=\sum_{1\leq n\leq \log^2 T} \frac{\Lambda(n)}{n^s}+O((\log T)^{2-2\sigma}),$$
where the implied constant depends at most on $\sigma_1$.
This completes the proof of the lemma.
\end{proof}

\section{Proof of Theorem~\ref{thm 1}}
We now proceed to prove Theorem~\ref{thm 1}, which we restate for the convenience of the reader.
\begin{thm 1}  
Assume RH. Let $\sigma_1$ and $\sigma_2$ be fixed with $1/2<\sigma_1<\sigma_2<1.$ 
 If $\rho'= \beta'+\di \gamma'$ is any critical point of the Riemann zeta function with
 $\sigma_1<\beta'<\sigma_2,$ then there is a positive constant $A$ depending on $\s_1$ and $\s_2$ such that
\be\notag
\log |\zeta(\rho')|\le A  \ {\frac{(\log \gamma')^{2-2\b'}}{\log \log \gamma'}}.
\ee
\end{thm 1}

\begin{proof}
By Lemma~\ref{log deriv}  we have
\begin{equation}\label{eq log deriv}
-\frac{\zeta'}{\zeta}(s)=\sum_{n\leq \log^2 T} \frac{\Lambda(n)}{n^s}+O((\log^2 T)^{1-\sigma}),   
\end{equation}
where $s=\sigma+\di t$, 
$1/2<\sigma_1\leq\sigma\leq\sigma_2<1$, and $4\leq T\leq t\leq2T.$
Integrating from $\beta'$ to $\infty$, where $\rho' =\beta'+i\gamma'$ is a critical point of the zeta function and $4\leq T \leq \gamma'\leq 2T$,
we see that  
$$\log\zeta(\rho')=\int_{\beta'}^{\infty} -\frac{\zeta'}{\zeta}(y+\di \gamma')dy
=\sum_{n\leq \log^2 T}
\frac{\Lambda(n)}{n^{\rho'} \log n} 
+O\bigg(\frac{(\log T)^{2-2\beta'}}{\log \log T}\bigg).$$
Taking  the real part of both sides, we obtain
$$\log|\zeta(\rho')|=\Re\sum_{ n\leq \log^2 T} \frac{\Lambda(n)}{n^{\rho'}\log n} +O\bigg(\frac{(\log T)^{2-2\beta'}}{\log \log T}\bigg).$$
From Lemma 3 we see that
$$\Re\sum_{2\leq n\leq \log^2 T} \frac{\Lambda(n)}{n^{\rho'} \log n}= O\bigg({\frac{(\log T)^{2-2\beta'}}{\log \log T}}\bigg),
$$
so that 
$$\log|\zeta(\rho')|=O\bigg(\frac{(\log T)^{2-2\beta'}}{\log \log T}\bigg).$$
Since $T\leq \gamma' \leq 2T$, it follows that
$$
\log |\zeta(\rho')|\ll{\frac{(\log \gamma')^{2-2\beta'}}{\log \log \gamma'}}
$$
where the implicit constant depends on $\s_1$ and $\s_2$.
\end{proof}	

 
\section{Lemmas for the proof of Theorem 2}   
 
The results in this section are all unconditional.

	\begin{lemma}\label{lemma1_2}
		Let $a(n)$ be a totally multiplicative function such that $|a(n)|\leq 1$ for all $n$. Then for all $x\geq 1$
	$$
		\sum_{n\leq x} a(n)\Lambda(n)n^{-s}=\sum_{p\leq x}\frac{a(p) \log p}{p^s-a(p)} +O\big(x^{1/2-\sigma}\big),
	$$
where $1/2<\sigma<1.$	
 \end{lemma}
 
		\begin{proof}
	We begin by observing that 
	\begin{align*}
	    \sum_{p\leq x}\frac{a(p) \log p}{p^{s}-a(p)}&=\sum_{p\leq x}\frac{a(p) \log p}{p^s}\bigg(\frac{1}{1-a(p)/p^s}\bigg)=\sum_{p\leq x}{p^{-s}}{a(p) \log p}\sum_{k=0}^{\infty}\frac{a(p^k)\log p}{p^{ks}}\\&=\sum_{p\leq x}\sum_{k=1}^{\infty}\frac{a(p^k)\log p}{p^{ks}}.
	\end{align*}
Thus, we see that
\be\label{a(p) sum 1}\sum_{n\leq x} a(n)\Lambda(n)n^{-s}=\sum_{p^k\leq x} 
a(p^k)(\log p) p^{-ks}=\sum_{p\leq x}\frac{a(p) \log p}{p^s-a(p)}-\sum_{k=2}^{\infty}\sum_{x^{1/k}<p\leq x}\frac{a(p^k)\log p}{p^{k\sigma}}.
\ee
Let $\theta(x)=\sum_{p\leq x}\log p.$ We know that $\theta(x)\ll x.$\\
Thus, for any $a>1$ and any $y_1,y_2$ such that $1<y_1<y_2$,
$$\sum_{y_1<p\leq y_2}\frac{\log p}{p^a}=\int_{y_1}^{y_2}\frac{d\theta(x)}{x^a}=O\big(y_1^{1-a}\big).$$
Using this with $x^{1/k}\geq 2$, that is, $2\leq k\leq \log x/\log 2$, we have	 
$$\sum_{x^{1/k}<p\leq x}\frac{\log p}{p^{k\sigma}}=O\big(x^{1/k-\sigma}\big).$$
And when $x^{1/k}<2$, that is, $k> \log x/\log 2$, we have
$$\sum_{2\leq p\leq x}\frac{\log p}{p^{k\sigma}}=O\big(2^{1-k\sigma}\big).$$
Combining these estimates, we find that the second term on the far  right in \eqref{a(p) sum 1} equals
\begin{align*}
 \sum_{k=2}^{\infty}\sum_{x^{1/k}<p\leq x}\frac{a(p^k)\log p}{p^{k\sigma}}&=\sum_{k=2}^{\log x/\log 2}\sum_{x^{1/k}<p\leq x}\frac{a(p^k)\log p}{p^{k\sigma}}+\sum_{k=\log x/\log 2}^{\infty}\sum_{2\leq p\leq x}\frac{a(p^k)\log p}{p^{k\sigma}}\\
&\ll  \sum_{k=2}^{\log x/\log 2} x^{1/k-\sigma} +\sum_{k=\log x/\log 2}^{\infty} 2^{1-k\sigma} \\
&\ll x^{1/2-\sigma} + x^{1/3-\sigma}\log x + x^{-\sigma} \ll x^{1/2-\sigma} .
\end{align*}
Inserting this on the right-hand side of \eqref{a(p) sum 1}, we obtain
$$
		\sum_{n\leq x} a(n)\Lambda(n)n^{-s}=\sum_{p\leq x}\frac{a(p) \log p}{p^s-a(p)} +O\big(x^{1/2-\sigma}\big).
$$
\end{proof} 

Next we  define some functions and parameters that will be helpful in constructing a root of $\zeta'(s)$ arbitrarily close to the line $\Re s=\sigma_1.$ 
Let $\sigma_1$ be fixed with $1/2<\sigma_1<1$ . 
Then for $x\geq 2, 0<c<1$ and $\s_1<\s<1$, we define auxiliary functions $V_x(s)$ and $W_x(s)$ as
\be \label{V}
V_x(s)=\displaystyle\sum_{n\leq x}{\Lambda(n)n^{-s}}
\ee
and
\be \label{W}
W_x(s)=\displaystyle\sum_{n\leq x}{b(n)\Lambda(n)n^{-s}},
\ee
where $b(n)$ is a totally multiplicative function such that $b(p)=1$ for all $p\leq cx$, and $b(p)=-1$ for all $p > cx$. 
Next we express $W_x(s)$ in terms of $V_x(s)$.
By Lemma \ref{lemma1_2},  if $\sigma_1<\sigma <1$, then
\begin{align*}
    W_x(s)&=\sum_{p\leq cx}\frac{\log p}{p^s-1}-\sum_{cx<p\leq x}\frac{\log p}{p^s+1}+O(x^{1/2-\sigma_1})\\
    &=2\sum_{p\leq cx}\frac{\log p}{p^s-1}-\sum_{p\leq cx}\frac{\log p}{p^s-1}-\sum_{cx<p\leq x}\frac{\log p}{p^s+1}+O(x^{1/2-\sigma_1})\\
    &=2\sum_{p\leq cx}\frac{\log p}{p^s-1}-\sum_{p\leq x}\frac{\log p}{p^s-1}+\sum_{cx <p\leq x}\frac{\log p}{p^s-1}-\sum_{cx<p\leq x}\frac{\log p}{p^s+1}+O(x^{1/2-\sigma_1}). 
    \end{align*}
 Applying  Lemma \ref{lemma1_2} again to the first two sums, and separately combining the last two sums, we find that   
 \be\label{Wx 1}
  \begin{split}  
    W_x(s) &=2V_{cx}(s)-V_x(s)+\sum_{cx<p\leq x}\frac{2\log p}{p^{2s}-1}+O(x^{1/2-\sigma_1})\\
    &=2V_{cx}(s)-V_x(s)+O(x^{1/2-\sigma_1}).
    \end{split}
    \ee
    
We now specify the value of $c$ in the definition of $W_x(s)$ as
\be\label{c}
  \log c= {{-\frac{\log 2}{1-\sigma_1}+\frac{\log{2}}{(1-\sigma_1)^2\log^ax}}}.
\ee 
We use this to show that $W_x(s)$ has a root near $\s_1$.
\begin{lemma}\label{lemma2_2}
Let $1/2<\s_1<1$,  let  $a >1$ be fixed, and let $c$ in the definition of $W_x(s)$ be given by \eqref{c}.
Then for  all  large $x$,
 $W_x(s)$ has a root at $$s=\sigma_1 +\frac{1}{\log^a x}+O\Big(\frac{1}{\log^{2a} x}\Big).$$
\end{lemma}

\begin{proof}
Writting $\psi(y)=\sum_{n\leq y}\Lambda(n)$
and applying the prime number theorem, we see  that 
\begin{equation}\label{eq1_2}
V_x(s)=\sum_{n\leq x}\frac{\Lambda(n)}{n^s}=\int_{1}^{x}\frac{d\psi(y)}{y^s}
=\frac{x^{1-s}}{1-s}+O\big(x^{1-\sigma_1}\exp(-c_1\sqrt{\log x})\big),
\end{equation}
where $c_1>0$ is an absolute constant.
Using this in \eqref{Wx 1}, we obtain
\be\label{W_x 1}
W_x(s)=\frac{2(cx)^{1-s}-x^{1-s} }{1-s}+O\big(x^{1-\sigma_1} \log^{-2a} x\big),
\ee
say.
Next we set  $s=\s_1+z$, where $|z| <2/\log^a x$  and $x$ is so large that 
$|1-s|\geq |1-\s_1|- 2/\log^a x >0$.

 Then 
\be\label{W_x 2}
\frac{W_x(s)}{x^{1-s}}=\frac{2 e^{(1-s) \log c}- 1}{1-s}+O ( \log^{-2a} x ).
\ee
Now 
\be\label{(1-s) log c}
\begin{split}
(1-s) \log c= &(1-\s_1-z) \Big({{-\frac{\log 2}{1-\sigma_1}+\frac{\log{2}}{(1-\sigma_1)^2\log^ax}}} \Big) \\
= &(1-\s_1)\Big(1-\frac{z}{1-\s_1} \Big) \Big( \frac{\log 2}{\sigma_1-1} \Big)\Big(1+ \frac{1}{(1-\s_1)\log^{a} x} \Big)  \\
=&-\log 2   \Big(1-\frac{z}{1-\s_1} + \frac{1}{(1-\s_1)\log^{a} x} + O(\log^{-2a} x) \Big).
\end{split}
\ee
Thus,
\be\label{W_x 3}
\begin{split}
\frac{W_x(s)}{x^{1-s}}=&\frac{ e^{ \log 2\big(\frac{z}{1-\s_1} - \frac{1}{(1-\s_1)\log^{a} x} 
+ O(\log^{-2a} x) \big)}- 1}{1-s}+O ( \log^{-2a} x )\\
=&\frac{  \log 2 \big(z  -  \log^{-a} x  \big)}{(1-\s_1)(1-s)}+O ( \log^{-2a} x ).
\end{split}
\ee
It follows that $W_x(s)$ has a root at
$$
s=\sigma_1 +\frac{1}{\log^a x}+O\Big(\frac{1}{\log^{2a} x}\Big).$$
\end{proof}

If  $p^{i\tau}$ is very near $1$ for $p \leq cx$ and $p^{i\tau}$ is very near $-1$ for $cx<p \leq x$, then $V_x (s + i\tau )$ will be close to $W_x (s)$. To show that such $\tau$ exist within a reasonable height, we need a sharp  form of Kronecker’s theorem concerning inhomogeneous Diophantine approximation. To this end, we follow closely the approach of Gonek and Montgomery in~\cite{paper}.


\begin{lemma}\label{lemma3_2}
Let $K$ be a positive integer and suppose $0<\delta\leq \frac{1}{2}$. There is a trigonometric polynomial $f(\theta)$ of the form
\be\label{f defn} 
f(\theta)=\sum_{k=0}^{K}c_ke(-k\theta)
\ee
such that $ max_\theta |f(\theta)|=f(0)=1$ and $|f(\theta)|\leq 2 e^{-\pi K \delta}$ for $\delta \leq \theta \leq 1-\delta.$
\end{lemma}	
\begin{proof}
This is Lemma  7 of   \cite{paper}.
\end{proof}

 The second moment
\begin{equation}\label{E:Defmu}
\mu = \int_0^1 |f(\theta)|^2\,d\theta
\end{equation}
appears below.  Since
\[
1 = |f(0)|^2 = \Big|\sum_{k=0}^K c_k\Big|^2 \le (K+1)\sum_{k=0}^K |c_k|^2
\]
by Cauchy's inequality, it follows that
\begin{equation}\label{E:muEst}
\frac1{K+1}\le \mu \le 1\,.
\end{equation}

For a given finite set $\mathscr P$ of primes $p$ and a given set of real numbers $\beta_p$ (considered modulo 1), we want to show that there exist  real numbers $t$ in prescribed intervals, such that $\big\|{t\log p/2\pi} -\beta_p \big\|<\delta$, where $\|x\|$ indicates the distance of $x$ from the nearest integer.	To accomplish this we define   
\begin{equation}\label{eqg(t)}
     g(t)=\prod_{p\in \mathscr P}\bigg|f\Big(\frac{t\log p}{2 \pi}-\beta_p\Big)\bigg|^2,	
\end{equation}
where $f$ is as in the previous lemma.

\begin{lemma}\label{lemma4_2}
	Let $\mathscr P$ be a set of primes not exceeding $x$. For each $p\in \mathscr P$ let a number 
	$\beta_p$ be given, and let 
	$K$,  $\mu$, and $ g(t)$ be as in  \eqref{f defn}, \eqref{E:Defmu},
 and \eqref{eqg(t)}. Then   for $Y$ any real number and $T\geq 4$,
	\be\label{int g 1}
		\int_{T}^{T+Y}g(t)dt=(Y+O(\exp(2Kx)))\mu^{card \mathscr P}.
	\ee
\end{lemma}	
\begin{proof}
The proof is almost identical to that of the proof of Lemma 8 in\cite{paper}.	
\end{proof}			
			
 From this point on, recalling that $a>1$ in some of the lemmas above, we let   
 \be\label{d_1, b}
 d_1=\sigma_1-1/2, \qquad  b>a+1>2,
\ee
and make the following choices for the parameters 
 $K$, $x$, and $\delta$:
\be\label{parameters}
  x=\frac{d_1\log T}{4(\log\log T)^b},\qquad K=\bigg[\frac{1}{2}\log^b x\bigg], \qquad \delta=\frac{1}{\log^{b-1} x}.
\ee

For each root $\rho=\beta + \di \gamma $ of the zeta function such that $\beta\geq (1+d_1)/2$,
 we remove from $[T,2T]$ those $\tau$ satisfying $|\gamma-\tau|\leq T^{d_1/4}+1$. We let $X$ denote the set of $\tau$ we have removed, and let $R$ denote the remaining set, so that 
 \be\label{R and X} 
 [T,2T]= R \cup X.
\ee
 By Theorem 9.19 (A) of Titchmarsh~\cite{Titch}, the number of roots with ordinates in $[T, 2T]$ is $\ll T^{1-d_1 /2}\log^{5}T.$
  Thus, the set  $X$ has  measure $\ll T^{1-d_1/4}\log^{5}T.$    
   Note that if $\tau \in R$ and $s=\sigma+it$ with   $|t|\leq 1$, then $min_\gamma |\gamma-t-\tau|\geq T^{d_1/4}.$
   In particular, this holds when $s$ is on or inside $\mathscr{C}_1$ and $\tau \in R$.

We next prove an analogue of Lemma~\ref{lemma4_2} for the integral $\int_{R}g(t)dt$.

\begin{lemma}\label{lemma5_2}
Under the same hypotheses as in Lemma~\ref{lemma4_2}, and with $R$ and $d_1$ as above, we have  
$$
\int_{R}g(t)dt=(T+O(T^{1-d_1/4}\log^{5}T) )\mu^{card \mathscr P}.
$$
\end{lemma}
\begin{proof}
We have 
\be\label{int g R}
\int_{R}g(t)dt=\int_{T}^{2T}g(t)dt -\int_{X}g(t)dt.
\ee
By Lemma~\ref{lemma4_2}
\be\notag  
\int_{T}^{2T}g(t)dt = (T+O(\exp(2Kx)))\mu^{card \mathscr P}.
\ee
 By our choice of the parameters $K$ and $x$ in \eqref{parameters}
we see that
 $2Kx\leq d_1(\log T)/4$, so that 
 \be\label{exp est} 
 \exp(2Kx)\leq T^{d_1/4}.
 \ee
Thus
\be\label{int g 2}
\int_{T}^{2T}g(t)dt = (T+O(T^{d_1/4}) )\mu^{card \mathscr P}.
\ee

 Next,  $X$ consists of $\ll T^{1-d_1 /2}\log^{5}T$ intervals, each  of length $\leq 2 T^{d_1/4}$. 
 Thus, by Lemma~\ref{lemma4_2} and \eqref{exp est}, each such interval contributes an amount
 $$
 \ll (T^{d_1/4}+\exp(2Kx)) \, \mu^{card \mathscr P}  \ll  T^{d_1/4} \mu^{card \mathscr P}.
 $$
 It follows that
\begin{equation}\notag
    \begin{split}
        \int_{X}g(t)dt&\ll   (T^{1-d_1/2}\log^{5}T )    \, (T^{d_1/4} \; \mu^{card \mathscr P} )\\
        &= (T^{1-d_1/4}\log^{5}T )\, \mu^{card \mathscr P} .  
    \end{split}
\end{equation}

Combining this and \eqref{int g 2}  in \eqref{int g R}, we obtain
 $$
 \int_{R}g(t)dt=(T+O(T^{1-d_1/4}\log^{5}T))\;\mu^{card \mathscr P}.
 $$
\end{proof}

The function $g(t)$ is large when the numbers $ \| t(\log p/2\pi)-\beta_p \|$ are small,  but  to obtain Kronecker’s theorem we need a peak function that is positive \emph{only} when all of these numbers are $<\delta$. To accomplish this we
define  
\be\label{h defn}
    h(t)=\displaystyle\prod_{p\leq x}\bigg|f\bigg(\frac{t\log p}{2 \pi}-\beta_p\bigg)\bigg|^2-\epsilon\sum_{p_1\leq x}\prod_{p\leq x,p\neq p_1}\bigg|f\bigg(\frac{t\log p}{2 \pi}-\beta_p\bigg)\bigg|^2,
\ee
where 
\be\label{epsilon} 
\epsilon=4e^{-2\pi K \delta}\ll x^{-3}.
\ee
It is easy to see that $h(t)>0$ only when 
$\big\|{t\log p/2\pi} -\beta_p \big\|<\delta$ for all $p\leq x$ (confer ~\cite{paper}).

\begin{lemma}\label{lemma6_2}
With $f(\theta)$ defined as in Lemma~\ref{lemma3_2}, $h(t)$ defined as above, and the choice of parameters in \eqref{parameters}, 
we have  for $T\geq 4$ and  any real $Y$, that
\be\label{h on Y inter}	
\int_{T}^{T+Y} h(t)dt= (Y+O(T^{d_1/4}))(1+O(x^{-1}) )\mu^{\pi(x)}  .
\ee
Moreover, for $R$ as in \eqref{R and X}, we have
\be\label{h on R}
\int_{R}h(t)dt = (1+O(x^{-1}) )\mu^{\pi(x)}T.
\ee
\end{lemma}
\begin{proof}
From the definition of $h(t)$, Lemma \ref{lemma4_2}, and \eqref{exp est}, we have
\be
\begin{split}
\int_{T}^{T+Y}h(t)dt
=&(Y+O(\exp(2Kx)))\mu^{\pi(x)}
+O \big(\epsilon \pi(x)(Y+O(\exp(2Kx)\mu^{\pi(x)-1}) )\big)\\
=&(Y+O(T^{d_1/4}) ) \mu^{\pi(x)}
+O \big(\epsilon \pi(x)(Y+O(T^{d_1/4}))\;\mu^{\pi(x)-1}\big).
\end{split}
\ee
Now $\pi(x)\ll{x/\log x}$,  $\epsilon\ll {x^{-3}}$ by \eqref{epsilon}, and by \eqref{E:muEst} and our choice of $K$ in \eqref{parameters},
${1/\mu}\ll (\log^b x) $.
It follows that  $\epsilon \pi(x)/\mu \ll \log^{b-1} x/x^2 \ll 1/x.$ This establishes \eqref{h on Y inter}.

Next we prove \eqref{h on R}. By \eqref{R and X},
\be\label{h int on R}
\int_{R}h(t)dt = \int_{T}^{2T} h(t)dt-\int_{X}h(t)dt.
\ee
By \eqref{h on Y inter} with $Y=T$,
\be\label{h on Y inter 2}	
\int_{T}^{2T} h(t)dt= (T+O(T^{d_1/4}))(1+O(x^{-1}) )\mu^{\pi(x)}  .
\ee
To estimate $\int_{X}h(t)dt$, recall that
$X$ consists of $\ll T^{1-d_1 /2}\log^{5}T$ intervals, each  of length $\leq 2 T^{d_1/4}$. 
Thus, by \eqref{h on Y inter}, each such interval contributes an amount
 $$
 \ll T^{d_1/4}(1+O(x^{-1})) \mu^{\pi(x)}   
 \ll  T^{d_1/4} \mu^{\pi(x)}.
 $$
 It follows that
\begin{equation}\notag
    \begin{split}
        \int_{X}h(t)dt&\ll   (T^{1-d_1/2}\log^{5}T )    \, (T^{d_1/4} \; \mu^{\pi(x)}\\
        &= (T^{1-d_1/4}\log^{5}T )\, \mu^{\pi(x)} .  
    \end{split}
\end{equation}
Combining this and \eqref{h on Y inter 2}  with \eqref{h int on R}, we obtain
 $$
 \int_{R}h(t)dt=T(1+O(x^{-1}) )\;\mu^{\pi(x)},
 $$
 which is \eqref{h on R}.
\end{proof}

These lemmas ensure that there are $t$ for  which the primes $p  \leq x$ behave as we want. However,  
the remaining primes $p >x$ could   make an unwanted contribution. The next lemma allows us to ensure this does not happen.

\begin{lemma}\label{lemma7_2}
Let $g(t)$ be as in \eqref{eqg(t)}, where $\mathscr P$ is the set of primes not exceeding $x$. 
For each $p>x$ let $b_p$ have the property that $|b_p| \leq 1/p^{\sigma_1}$.
Then
$$\int_{T}^{2T}g(t)\bigg|\sum_{x <p\leq T^{d_1/4}}\frac{b_p}{p^{\di t}}\bigg|^2dt\ll T\mu^{\pi(x)}\frac{x^{1-2\sigma_1}}{\log x},$$ 
where the implied constant depends on $\sigma_1$.
The same bound holds a fortiori for the integral over $R$.	
\end{lemma}

\begin{proof}
By \eqref{f defn} we see that
\[
\prod_{p\in\mathscr P} f\Big(t\frac{\log p}{2\pi} - \beta_p\Big)
= \prod_{p\in\mathscr P}\Big(\sum_{k=0}^Kc_ke(k\beta_p)p^{-ikt}\Big)
= \sum_{n\in \mathscr N} a_nn^{-it}
\]
where $\mathscr N$ is the set of positive integers composed entirely of primes
in $\mathscr P$, with multiplicities not exceeding $K$ and
\[
a_n = \prod_{\substack{p\in\mathscr P\\ p^k\|n}}c_ke(k\beta_p)\,.
\]
Here the product is extended over all members of $\mathscr P$, not just those
dividing $n$.
We note that a positive integer $m$ has at most one decomposition $m = np$
with $n\in \mathscr N$ and $p > x$.  
Let the numbers $C_m$ be determined by the identity
\[
\Big(\sum_{n\in\mathscr N} a_nn^{-it}\Big)\Big(\sum_{x<p \leq T^{d_1/4}}b_pp^{-it}\Big)
= \sum_m C_m m^{-it}\,.
\]
 Montgomery \& Vaughan \cite{paper2} have shown that
if $\sum_m |C_m|<\infty$, then
\[
\int_T^{2T}\Big|\sum_{m=1}^\infty C_m m^{-it}\Big|^2\,dt
= \sum_{m=1}^\infty |C_m|^2(T+O(m))\,.
\]
In the main term we have
\[
\sum_{m=1}^\infty |C_m|^2 = \Big(\sum_{n\in \mathscr N} |a_n|^2\Big)
\Big(\sum_{x<p \leq T^{d_1/4}}|b_p|^2\Big)\,.
\]
The sum over $n$ is $\mu^{\pi(x)}$, and the sum over $p$ is
$\ll \sum_{x<p \leq T^{d_1/4}} p^{-2\sigma_1} \ll x^{1-2\sigma_1}/\log x$.  In the error term we have
\[
\sum_{m=1}^\infty m|C_m|^2 = \Big(\sum_{n\in\mathscr N}n|a_n|^2\Big)
\Big(\sum_{x<p \leq T^{d_1/4}}p|b_p|^2\Big)\,.
\]
For $n\in\mathscr N$ we have $n\le \exp(2Kx) \leq T^{d_1/4}$, so the sum over $n$
here is $\ll \mu^{\pi(x)}T^{d_1/4}$.  
The sum over $p$ is
\[
\le \sum_{x<p \leq T^{d_1/4}}p^{1-2\sigma_1} \leq T^{d_1(1-\sigma_1)/2}/\log T.
\] 
Combining our estimates, we obtain
\be
\begin{split}
\int_{T}^{2T}g(t)\bigg|\sum_{x <p\leq T^{d_1/4}}\frac{b_p}{p^{\di t}}\bigg|^2dt 
\ll & T\mu^{\pi(x)}\frac{x^{1-2\sigma_1}}{\log x}+T^{d_1/4}\mu^{\pi(x)}\frac{T^{d_1(1-\sigma_1)/2}}{\log T}\\
\ll &T\mu^{\pi(x)}\frac{x^{1-2\sigma_1}}{\log x} 
\end{split}
\ee
by our choice of $x$.
This completes the proof of the lemma.
\end{proof}

 \begin{lemma}\label{lem 12}
Let $W_x(s)$ be as in \eqref{W} with $x=(d_1 \log T)/4(\log \log T)^b,  d_1 =\sigma_1-1/2$, and $T\geq 4$.
Then for $s=\sigma+it$ with $\sigma_1 \leq \sigma$ and $|t|\leq 1$,  and for $\tau \in R$, where $R$ is defined just before \eqref{R and X}, we have 
 \begin{equation}\label{eq15_2}
\begin{split}
-\frac{\zeta'}{\zeta}(s+\di\tau)
=\sum_{n\leq T^{d_1/4}}w(n)\Lambda(n)n^{-s-\di\tau}+O(T^{-d_1^2/16}),    
\end{split}
\end{equation}
where  
\be\label{w}
 w(u)=
\begin{cases}
	1  &\qquad \hbox{if}\qquad 1\leq u\leq T^{d_1/8},\\
	1  -\frac{\log(u/y)}{\log y} &\qquad  \hbox{if}\qquad  T^{d_1/8}<u\leq T^{d_1/4},\\
	 0 &\qquad  \hbox{if}\qquad  u>T^{d_1/4}.				
\end{cases}
\ee	        

 \end{lemma}  
 \begin{proof}
By \eqref{MV expl frml} with $x=y=T^{d_1/8}, T\geq 4, d_1 =\sigma_1-1/2, \sigma\geq \sigma_1$, and $\tau \varepsilon R$, we have
\begin{equation}\label{eqmv}
    \begin{split}
     -\frac{\zeta'}{\zeta}(s+\di\tau)&=-\sum_{n\leq y^2} w(n)\frac{\Lambda(n)}{n^{s+\di\tau}}+\frac{y^{2(1-s+\di\tau)}-y^{1-s-\di\tau}}{(1-s-\di\tau)^2\log y}-\sum_{\rho}\frac{y^{2(\rho-s-\di\tau)}-y^{\rho-s-\di\tau}}{(\rho-s-\di\tau)^2\log y}\\&-\sum_{k=1}^{\infty}\frac{y^{2(-2k-s-\di\tau)}-y^{-2k-s-\di\tau}}{(2k+s+\di\tau)^2\log y}.
\end{split}
\end{equation}	
%

Since $R\subseteq [T, 2T]$, the second term on the right-hand side of \eqref{eqmv} becomes
$$
\frac{y^{2(1-s-\di\tau)}-y^{1-s-\di\tau}}{(1-s-\di\tau)^2\log y}\ll \frac{T^{d_1(1-\sigma_1)/4}}{d_1 T^2\log T}
\ll  T^{-d_1/8-d_1^2/4-2}.
$$
We  split the third term on the right-hand side of \eqref{eqmv} into two sums $P$ and $Q$, where
$P$ is over the zeros with $\beta\geq (1+d_1)/2$, and $Q$ is over the zeros with $\beta< (1+d_1)/2$.
For Q we have
\begin{align*}
Q=  \sum_{\beta<(1+d_1)/2}\frac{y^{2(\rho-s-\di\tau)}-y^{\rho-s-\di\tau}}{(\rho-s-\di\tau)^2\log y}&\ll\sum_{\beta<(1+d_1)/2}\frac{y^{2(\beta-\sigma)}+y^{\beta-\sigma}}{[(\beta-\sigma)^2+(\gamma-t-\tau)^2]\log y}\\
&\ll \frac{y^{-d_1/2} }{\log y} \sum_{\gamma}\frac{1}{d_1^2+(\gamma-(t+\tau))^2}\\
&\ll  \frac{y^{-d_1/2} \log T}{d_1^2\log y} \\
&\ll y^{-d_1/2}=T^{-d_1^2/16} . 
\end{align*}
We also have
\begin{align*}
P=  \sum_{\beta\geq (1+d_1)/2}\frac{y^{2(\rho-s-\di\tau)}-y^{\rho-s-\di\tau}}{(\rho-s-\di\tau)^2\log y}
&\ll\sum_{\beta\geq(1+d_1)/2}\frac{y^{2(1-\sigma)}+y^{1-\sigma}}{(\gamma-t-\tau)^2\log y}\\
&\ll \frac{y^{1-2d_1}\log T}{T^{d_1/4}\log y} 
\ll \frac{y^{1-2d_1}}{T^{d_1/4} } =T^{-d_1/8-d_1^2/4}.  
\end{align*}
The last term on the right-hand side of \eqref{eqmv}
is
$$\sum_{k=1}^{\infty}\frac{y^{2(-2k-s-\di\tau)}-y^{-2k-s-\di\tau}}{(2k+s+\di\tau)^2\log y}
\ll \frac{y^{-2-\sigma_1}}{T\log y}\ll T^{-1-5d_1/16 -d_1^2/8}.$$
Combining all these estimates, we find that for $\sigma\geq \sigma_1 $ and $\tau \in R$,
\begin{equation}\label{eq13_2}
    -\frac{\zeta'}{\zeta}(s+\di\tau) 
    =-\sum_{n\leq T^{d_1/4}} w(n)\frac{\Lambda(n)}{n^{s+\di\tau}}+O(T^{-d_1^2/16}).
\end{equation}
\end{proof}


\section{The proof of Theorem~\ref{thm 2}}	

Let $V_x(s)$ and $W_x(s)$ be  as in \eqref{V} and \eqref{W}.
Let 
$$
\mathscr{C}_0 =\Big\{ s=\sigma_1+\frac{1+ e^{\di \theta}/2 }{\log^a x} \; : \;  0\leq \theta \leq 2\pi \Big\}
$$
and
$$\mathscr{C}_1 =\Big\{ s=\sigma_1+\frac{1+ e^{\di \theta}/4 }{\log^a x} \; : \;  0\leq \theta \leq 2\pi \Big\},
$$
where, as previously, $a>1$.
Also, let $c$ be as in \eqref{c} and $\delta$ as in \eqref{parameters}.
Suppose that  $\tau$ is a real number such that 
$$
\Big\|{\tau \frac{\log p}{2\pi}} \Big\|<\delta  \quad\hbox{for} \quad  p\leq cx,
$$
$$
\Big\|{\tau \log p/2\pi}+{1/2} \Big\|<\delta   \quad\hbox{for} \quad cx <p\leq x,
$$ 
and let $\mathscr{G} $ be the set of $\tau$ such that both inequalities hold.

By Lemma~\ref{lemma1_2}, \eqref{V}, and \eqref{W} we see that for $\sigma>\sigma_1$

\be
V_x(s) =\sum_{p\leq x} \frac{\log p}{p^s -1}  +O(x^{1/2-\sigma_1})
  \ee
  and
\be
W_x(s) =\sum_{p\leq x} \frac{\log p}{b(p)p^s -1}  +O(x^{1/2-\sigma_1}).
  \ee
Hence
\be
\begin{split}
V_x(s+\di\tau)-W_x(s) = &
\sum_{p\leq x} \log p \Big(\frac{1}{p^{s+i\tau} -1}  - \frac{1}{b(p)p^s -1} \Big)
+O(x^{1/2-\sigma_1}) \\
\ll&\sum_{p\leq x} \frac{\log p}{p^{\sigma_1}} |p^{i\tau} - b(p)|
+O(x^{1/2-\sigma_1}).
\end{split}
\ee
Now note that for real $\theta$, if $\| \theta \|$ is the distance between $\theta$ and the nearest integer, then
$$
|e^{2\pi i \theta} -1| =2|\sin \pi \theta | \leq 2\pi \|\theta\| .
$$
Thus, taking 
$$
\theta=\theta_p=
\begin{cases}
\displaystyle \frac{\tau \log p}{2\pi} & \quad\hbox{if} \quad p\leq cx,\\
\displaystyle  \frac{\tau \log p}{2\pi}+\frac12 & \quad\hbox{if} \quad cx<p\leq x,
\end{cases}
$$
we see that  if $\tau\in \mathscr G$, then $|p^{i\tau} -b(p)|\leq 2\pi \delta$ for every $p\leq x$.
Thus, for $\sigma\geq \sigma_1$,
\be\label{V-W}
|V_x(s+\di\tau)-W_x(s)|\ll 
 \delta \sum_{p\leq x} \frac{\log p}{p^{\sigma_1}} 
+O(x^{1/2-\sigma_1})
\ll 
\frac{x^{1-\sigma_1}}{\log ^{b-1}x},
\ee
since $\delta=1/(\log x)^{b-1}$.
 
%
 
Define $$L_x(s)=\sum_{1<n\leq x}\frac{b(n)\Lambda(n)}{n^s \log n}.$$
Then we have
\begin{align}\notag
 \bigg|\sum_{1<n\leq x}\frac{\Lambda(n)}{\log n}n^{-s-\di\tau} -L_x(s)\bigg|
\leq \sum_{1<n\leq x} \frac{\Lambda(n)}{n^\sigma \log n}  | n^{-i\tau} -b(n)|   .
\end{align}
If $n=p^k$, then $ | n^{-i\tau} -b(n)| = |p^{ik\tau}-b(p)^k|\leq k | p^{i\tau} -b(p)|$. Thus, for $\sigma\geq \sigma_1$
and $\tau \in \mathscr G$, the expression here is
$$
\leq 2\pi \delta \sum_{p^k\leq x}\frac{1}{p^{k\sigma_1}} \ll \delta \frac{x^{1-\sigma_1}}{\log x}
\ll \frac{x^{1-\sigma_1}}{\log^b x}
$$
by our choice of $\delta$ in \eqref{parameters}.
That is, for $\sigma\geq \sigma_1$ and $\tau \in \mathscr G$,
\begin{align}\label{sum L_x}
 \bigg|\sum_{1<n\leq x}\frac{\Lambda(n)}{\log n}n^{-s-\di\tau} -L_x(s)\bigg|
\ll \frac{x^{1-\sigma_1}}{\log^b x}.\end{align}

Furthermore,  when $\Re s\geq \sigma_1$,
\begin{align*}
L'_x(s)=-\sum_{n\leq x}\frac{b(n)\Lambda(n)}{n^s}
\ll \sum_{n\leq x}\frac{\Lambda(n)}{n^{\sigma_1}}
\ll x^{1-\sigma_1}.
\end{align*}
Thus, if $s$ is on or inside  $\mathscr{C}_1$, then
\be\label{L(s)}
L_x(s)=L_x(\sigma_1)+O\bigg(\frac{x^{1-\sigma_1}}{\log^ a x}\bigg).
\ee
At $\sigma_1$,
\begin{equation}\label{eq6_2}
L_x(\sigma_1)=\sum_{1<n\leq x}\frac{b(n)\Lambda(n)}{n^{\sigma_1} \log n}
=2\sum_{p\leq cx}p^{-\sigma_1}-\sum_{p\leq x}p^{-\sigma_1}+O(1).    
\end{equation}
By the prime number theorem, for $u\geq 2$,
 $$\pi(u)=\frac{u}{\log u}+\frac{u}{\log^2 u}+O\Big(\frac{u}{\log^3 u}\Big).$$
From this and integration by parts, we see that for $y\geq 2$,
\begin{equation}\label{eq7_2}
   \sum_{p\leq y}p^{-\sigma_1}=\int_{2^{-}}^{y}\frac{d\pi(u)}{u^{\sigma_1}}=\frac{y^{1-\sigma_1}}{(1-\sigma_1)\log y}+\frac{y^{1-\sigma_1}}{(1-\sigma_1)^2\log^2 y}+O\bigg(\frac{y^{1-\sigma_1}}{\log^3 y}\bigg). 
\end{equation}
Using this in  \eqref{eq6_2}, we obtain
$$L_x(\sigma_1)=\frac{x^{1-\sigma_1}}{1-\sigma_1}\bigg(\frac{2c^{1-\sigma_1}}{\log cx}-\frac{1}{\log x}\bigg)+\frac{x^{1-\sigma_1}}{(1-\sigma_1)^2}\bigg(\frac{2c^{1-\sigma_1}}{\log^2 cx}-\frac{1}{\log^2 x}\bigg)+O\bigg(\frac{x^{1-\sigma_1}}{\log^3x}\bigg).$$
Now, from \eqref{(1-s) log c} with $s=\sigma_1$, that is, $z=0$, we have 
\be
 2c^{1-\sigma_1} = 1 - \frac{\log 2}{(1-\s_1)\log^{a} x} + O(\log^{-2a} x) =1+O({\log^{-a}x}).
\ee
Moreover,
$$
\frac{1}{\log cx} =\frac{1}{\log x}\Big(1 -\frac{\log c}{\log x}+ O\Big(\frac{1}{\log^2 x}\Big)\Big) .
$$
Hence, 
\begin{align*}
\frac{2c^{1-\sigma_1}}{\log cx}-\frac{1}{\log x}
=& \frac{(1+O(\log^{-a}x))}{\log x}\Big(1 -\frac{\log c}{\log x}+ O\Big(\frac{1}{\log^2 x}\Big)\Big) - \frac{1}{\log x}\\
=&-\frac{\log c}{\log^2x}+O\Big(\frac{1}{\log^{a+1} x}\Big) +O\Big(\frac{1}{\log^3 x}\Big)\\
=&-\frac{\log c}{\log^2x}+O\big( \log^{-\min(a+1, 3)} x \big)
\end{align*}
Similarly, one sees that  
$$\frac{2c^{1-\sigma_1}}{\log^2 cx}-\frac{1}{\log^2 x}=O(\log^{-3}x).$$
Thus, we find that 
$$
L_x(\sigma_1)=\frac{-x^{1-\sigma_1}\log c}{(1-\sigma_1)\log^2 x}
+O\bigg( \frac{x^{1-\sigma_1}}{ \log^{\min(a+1, 3)} x } \bigg).
$$ 
It now follows from \eqref{L(s)} that for $s$ on or inside  $\mathscr C_1$,
$$
L_x(s)=\frac{-x^{1-\sigma_1}\log c}{(1-\sigma_1)\log^2 x}
+O\bigg( \frac{x^{1-\sigma_1}}{ \log^{\min(a, 3)} x } \bigg).
$$
From this and \eqref{sum L_x}  we see that for $s$ on or inside  $\mathscr C_1$, and $\tau\in \mathscr G$,
 \begin{align}\label{sum L_x}\notag
\Re \sum_{1<n\leq x}\frac{\Lambda(n)}{\log n}n^{-s-\di\tau} =&\frac{-x^{1-\sigma_1}\log c}{(1-\sigma_1)\log^2 x}
+O\Big( \frac{x^{1-\sigma_1}}{ \log^{\min(a, 3)} x } \Big)
+O\Big(\frac{x^{1-\sigma_1}}{\log^b x}\Big)\\
=&\frac{-x^{1-\sigma_1}\log c}{(1-\sigma_1)\log^2 x}
+O\Big( \frac{x^{1-\sigma_1}}{ \log^{\min(a, 3)} x } \Big).
\end{align}
since $b>a+1$ and for the above inequality to hold true we need $a>2$.

Let
$$
T_x(s) = \sum_{x<n\leq T^{d_1/4}} \frac{w(n)\Lambda(n)}{\log n} n^{-s},
$$
where $\displaystyle x= \frac{d_1\log T}{4(\log\log T)^b}$.
Then by Lemma~\ref{lem 12}, for $s=\sigma+it$, with $\sigma_1 \leq \sigma$ and $|t|\leq 1$, and for $\tau \in R$,  we have 
\begin{equation}\label{ }
\log {\zeta}(s+\di\tau)
=\sum_{n\leq x} \frac{\Lambda(n)}{\log n}n^{-s-\di\tau}+T_x(s+i\tau) +O(T^{-d_1^2/16}),    
\end{equation}
and
\begin{equation}\label{ }
-\frac{\zeta'}{\zeta}(s+\di\tau)
=\sum_{n\leq x}\Lambda(n)n^{-s-\di\tau} -T_{x}^{'}(s+i\tau) +O(T^{-d_1^2/16}).    
\end{equation}
Now suppose that  $\tau \in R \cap \mathscr G$, and that
 \be\label{T_x}
 T_{x}(s+i\tau) \ll \frac{1}{\log x} , \qquad T_{x}^{'}(s+i\tau) \ll 1 
 \ee
for $s \in \mathscr{C}_1$. Then by \eqref{sum L_x} and \eqref{V-W}.
\begin{equation}\label{Re log zeta}
\Re \log {\zeta}(s+\di\tau)
=\frac{-x^{1-\sigma_1}\log c}{(1-\sigma_1)\log^2 x}
+O\bigg( \frac{x^{1-\sigma_1}}{ \log^{\min(a, 3)} x } \bigg)
\end{equation}
and
\begin{equation}\label{zeta'/zeta}
-\frac{\zeta'}{\zeta}(s+\di\tau)
=W_x(s) + O\Big(\frac{x^{1-\sigma_1}}{\log^{b-1} x}\Big)
\end{equation}
for $s$  on or inside $\mathscr{C}_1$.

 Recall from \eqref{W_x 3} that if $s=\sigma_1+z$ and $|z|< 2/\log^a x$, then
\be 
\begin{split}
\frac{W_x(s)}{x^{1-s}}
=&\frac{  \log 2 }{(1-\s_1)(1-s)}(z  -  \log^{-a} x  )+O ( \log^{-2a} x ).
\end{split}
\ee
Thus, for $s$  on or inside $\mathscr{C}_1$, 
\be\label{W_x 4} 
\begin{split}
 W_x(s) 
= x^{1-\sigma_1} \Big(\frac{e^{i\theta} \log 2}{4(1-\s_1)^2 \log^{a} x  }+O \Big( \frac{1}{\log^{2a} x }\Big)\Big)
\end{split}
\ee
(changed the error term to 2a)\\
Since $a>1$, we see from this that  the argument of $W_x(s)$ increases by $2\pi$ as $s$ traverses $\mathscr C_1$.
 Thus, by \eqref{zeta'/zeta}, $\zeta'(s)$ has a zero $\rho'$ in $\mathscr{C}_1+i\tau$, and from \eqref{Re log zeta}
 we see that
 \be
 \log | {\zeta}(\rho')|
\geq \Big(\frac{-\log c}{1-\sigma_1}+o(1)\Big)\frac{ x^{1-\sigma_1}}{\log^2 x} 
\ee
Now by \eqref{c},
$$ \log c = \log 2/(\sigma_1-1) (1+o(1))$$
so, using this and substituting 
$ {d_1\log T}/{4(\log\log T)^b}$ for $x$, we obtain
\begin{equation}\label{eq14_2}
\log|\zeta(\rho')| \geq  (1+o(1) ) \Big(\frac{(\sigma_1-1/2)^{1-\s_1}\log 2}{(1-\s_1)^2 4^{1-\s_1}}\Big)
\frac{( \log T)^{1-\sigma_1} }{(\log \log T)^{ 2+b(1-\sigma_1) } }.
\end{equation}
Since $T\leq \gamma'=\Im \rho' \leq 2T$,  \eqref{eq14_2} also holds with $T$ replaced by $\g'$, which is
\eqref{eqn thm 2}.

Thus we get the following equation.
$$\log|\zeta(\rho')|\geq \big((B(\s_1)-\epsilon)\big)\frac{(\log\gamma')^{1-\beta'}}{(\log\log\gamma')^{2+b(1-\sigma_1)}},$$ 
where
$$
  B(\sigma_1)  =\frac{(\sigma_1-1/2)^{1-\s_1}\log 2}{(1-\s_1)^2 4^{1-\s_1}}.
  $$

To complete the proof of Theorem~\ref{thm 2} it  only  remains to show the existence of a $\tau \in R\cap\mathscr G$ satisfying \eqref{T_x}. 
Recall that
$$T_x(s+\di\tau)=\sum_{x<p\leq T^{d_1/4}}w(p) p^{-s-\di\tau},$$ and
$$T'_x(s+\di\tau)=-\sum_{x<p\leq T^{d_1/4}}w(p)(\log p) p^{-s-\di\tau}.$$
We will prove  that there is a  constant $C_1$ that is independent of $x$ such that on the circle $\mathscr{C}_0 $, we have
$$\oint\limits_{\mathscr{C}_0}
|T_x(z+\di\tau)^2|dz|\leq \frac{C_1}{\log^9 x }.$$
Then, if $s$ is on or inside $\mathscr{C}_1$, we see from Cauchy's formula and the Cauchy-Schwarz inequality that
$$T_x(s+\di\tau)=\frac{1}{2\pi\di}\oint\limits_{\mathscr{C}_0} \frac{T_x(z+\di\tau)}{z-s}dz
\ll (\log x)^{a/2}\sqrt{\oint\limits_{\mathscr{C}_0} |T_x(z+\di\tau)^2|dz|}\ll \frac{1}{\log^{(9-a)/2} x },$$
and
$$T'_x(s+\di\tau)=\frac{1}{2\pi\di}\oint\limits_{\mathscr{C}_0} \frac{T_x(z+\di\tau)}{(z-s)^2}dz\ll (\log x)^{3a/2}\sqrt{\oint\limits_{\mathscr{C}_0} |T_x(z+\di\tau)^2|dz|} \ll \frac{1}{\log^{(9-3a)/2} x }.$$
From these estimates  we conclude that
$$\sum_{x< n\leq T^{d_1/4}}\frac{w(n)\Lambda(n)}{n^{s+\di\tau}\log n}=\sum_{x< p\leq T^{d_1/4}}\frac{w(p)}{p^{s+\di\tau}}+O\big(x^{1/2-\sigma_1}\big)=T_x(s+\di\tau)+O(x^{1/2-\sigma_1})\ll\frac{1}{\log^{(9-a)/2} x},$$\par
and $$\sum_{x< n\leq T^{d_1/4}}\frac{w(n)\Lambda(n)}{n^{s+\di\tau}}=\sum_{x< p\leq T^{d_1/4}}\frac{w(p)\log p}{p^{s+\di\tau}}+O\big(x^{1/2-\sigma_1}\big)=-T'(s+\di\tau)+O(x^{1/2-\sigma_1})\ll\frac{1}{\log^{(9-3a)/2}}.$$
In other words, the inequalities in \eqref{T_x} are true.

Now define $$h^+(t)=max{\{0,h(t)\}}.$$
Since $h(t)\leq h^{+}(t)\leq g(t)$, using Lemma \ref{lemma6_2}, we know that
$$\int_{R}h^{+}(t)dt=\mu^{\pi(x)}T\big(1+O(1/x)\big).$$

On the other hand we have
\begin{align*}
\int_{R}h^{+}(t)\Bigg(  \oint\limits_{\mathscr{C}_0} |T_x(z+\di t)|^2|dz|      \Bigg)dt&\leq \int_{R}g(t)\Bigg(  \oint\limits_{\mathscr{C}_0} |T_x(z+\di t)|^2|dz|      \Bigg)dt\\
&=\oint\limits_{\mathscr{C}_0}\int_{R}g(t)|T_x(z+\di t)|^2 dt |dz|\\
&\ll\mu^{\pi(x)}T\frac{x^{1-2\sigma_1}}{\log x} \\
&\ll\mu^{\pi(x)}T\frac{x^{1-2\sigma_1}}{\log^4 x}.\\
\text{Thus we conclude that}
\oint\limits_{\mathscr{C}_0} |T_x(z+\di t)|^2|dz|\leq \frac{\tilde C}{\log^9 x}.
\end{align*}
Choosing $b=3+\epsilon''$, we get
$$
\log|\zeta(\rho')| \geq  ((B(\s_1)-\epsilon)\big)
\frac{( \log T)^{1-\sigma_1} }{(\log \log T)^{ 5-3\beta+\epsilon' } },
$$
where $\epsilon'=\epsilon''(1-\sigma).$

Thus completing the proof.

\newpage

\section{Lemmas for the proof of Theorem 3}  
 
We assume RH for the proof of Theorem 3.
\begin{lemma}\label{lemma1}
Let $\sigma_1$ and $\sigma_2$ be such that $1/2<\sigma_1<\sigma_2<1$, then we have 
  $$\Re\sum_{1<n\leq x}\frac{\Lambda(n)}{n^s\log n}\geq -\frac{x^{1-\sigma}}{(1-\sigma)\log x}+O\Bigg(\frac{x^{1-\sigma}}{\log^2 x}\Bigg).$$ where $\sigma_1<\sigma<\sigma_2$ and the big-oh constant depends on $\sigma_1$ and $\sigma_2$.
\begin{proof}
We have the following inequality by Lemma 1 and Lemma 2 in part 1.
\begin{align*}
\Re\sum_{2\leq n\leq x} \frac{\Lambda(n)}{n^s \log n}&\geq -\sum_{2\leq n\leq x} \frac{\Lambda(n)}{n^\sigma \log n}=-\frac{x^{1-\sigma}}{(1-\sigma)\log x}+O\bigg({\frac{x^{1-\sigma}}{{\log^2 x}}}\bigg).
\end{align*}

\end{proof}
\end{lemma}
\section{Proof of theorem 3}
\begin{theorem}
Assume RH. Let $\sigma_1$ and $\sigma_2$ be fixed with $1/2<\sigma_1<\sigma_2<1.$ 
 If $\rho'= \beta'+\di \gamma'$ is any critical point of the Riemann zeta function with
 $\sigma_1<\beta'<\sigma_2,$ then there is a positive constant $C$ depending on $\s_1$ and $\s_2$ such that
 $$\log|\zeta(\rho')|\geq\frac{-C(\log \gamma')^{2-2\beta'}}{\log \log \gamma'}.
$$
\end{theorem}
\begin{proof}
By Lemma 4 in part 1 we have
\begin{equation}\label{eq log deriv}
-\frac{\zeta'}{\zeta}(s)=\sum_{1\leq n\leq \log^2 T} \frac{\Lambda(n)}{n^s}+O\big((\log^2 T)^{1-\sigma}\big)   
\end{equation}
for $\sigma_1\leq\sigma\leq\sigma_2$ and  let $4\leq T\leq t\leq2T.$
Integrating from $\beta'$ to $\infty$, where $\rho' =\beta'+i\gamma'$ is a critical point of the zeta function and $4\leq T \leq \gamma'\leq 2T$ 
we see that 
$$\log\zeta(\rho')=\int_{\beta'}^{\infty} -\frac{\zeta'}{\zeta}(y+\di \gamma')dy
=\sum_{2\leq n\leq \log^2 T}
\frac{\Lambda(n)}{n^{\rho'} \log n} 
+O\bigg(\frac{(\log T)^{2-2\beta'}}{\log \log T}\bigg).$$
Taking  the real part of both sides, we obtain
\begin{equation}\label{logzeta3}
    \log|\zeta(\rho')|=\Re\sum_{2\leq n\leq \log^2 T} \frac{\Lambda(n)}{n^{\rho'}\log n} +O\bigg(\frac{(\log T)^{2-2\beta'}}{\log \log T}\bigg).
\end{equation}
Because of the inequality above we can use Lemma 3 from part 1.
We can say that
\begin{equation}\label{approximation3}
    \Re\sum_{2\leq n\leq \log^2 T} \frac{\Lambda(n)}{n^{\rho'} \log(n)}\geq-\frac{(\log T)^{2-2\beta'}}{2(1-\beta')\log \log T} +O\bigg({\frac{(\log T)^{2-2\beta'}}{(1-\beta')^2 (\log \log T)^2}}\bigg).
\end{equation}

Combining \eqref{logzeta3} and \eqref{approximation3}, we get
$$\log|\zeta(\rho')|\geq-\frac{C(\log T)^{2-2\beta'}}{\log \log T} $$ for some constant $C.$
Since $4\leq T \leq \gamma'\leq 2T$ , substituting $\gamma'$ in the above equation we get
\\$$\log|\zeta(\rho')|\geq\frac{-C(\log \gamma')^{2-2\beta'}}{\log \log \gamma'} $$for some constant $B$
that depends on $\sigma_1$ and $\sigma_2$.\\
\end{proof}	

\newpage 
\section{Lemmas for the proof of Theorem 4}
We do not assume RH for the proof of Theorem 4.\\
Let $\sigma_1$ be such that $1/2<\sigma_1<1$ and $\sigma_1<\sigma<1.$
For $x\geq 1$, define $$Z_x(s)=\sum_{n\leq x}\frac{c(n)\Lambda(n)}{n^s},$$
where $c(n)$ is a totally multiplicative function such that $c(p)=-1$ for $p\leq cx$ and $c(p)=1$ for $p>cx$, where $$\log c={{-\frac{\log 2}{1-\sigma_1}+\frac{\log{2}}{(1-\sigma_1)^2\log^ax}}}.$$ 
Here we want $\sigma$ to be very close to $\sigma_1$ for large values of $x.$
\begin{lemma}\label{lemma11_4}
Let $\sigma_1$ be such that $1/2<\sigma_1<1$,  and $a$ be fixed. Then
for all large values of $x$, $Z_x(s)$ has a root at $$s=\sigma_1 +\frac{1}{\log ^ax}+O\bigg(\frac{1}{\log^{2a} x}\bigg)$$ where the big-oh constant depends on $\sigma_1.$
\begin{proof}
Let $x\geq 2$.
We set $z=s-\sigma_1$  where $|z| <2/\log^a x$  and $x$ is so large that 
$|1-s|\geq |1-\s_1|- 2/\log^a x >0$.
\\
By choosing sufficiently large $x$ we make sure that $s$ is arbitrarily close to $\sigma_1$.
Using Lemma \ref{lemma1_2} in part 2, we get 
\begin{align*}
Z_x(s)&=-\sum_{p\leq cx}\frac{\log p}{p^s+1}+ \sum_{cx<p\leq x}\frac{\log p}{p^s-1}+ O(x^{1/2-\sigma_1})\\
&=-\sum_{p\leq cx}\frac{\log p}{p^s+1}-\sum_{ p\leq cx}\frac{\log p}{p^s-1}+ \sum_{ p\leq x}\frac{\log p}{p^s-1}+ O(x^{1/2-\sigma_1})\\
&=-\sum_{p\leq cx}\frac{\log p}{p^s+1}+\sum_{ p\leq cx}\frac{\log p}{p^s-1}-2\sum_{p\leq cx}\frac{\log p}{p^s-1}+ \sum_{ p\leq x}\frac{\log p}{p^s-1}+ O(x^{1/2-\sigma_1})\\
\end{align*}
 Applying  Lemma \ref{lemma1_2} again to the first two sums, and separately combining the last two sums, we find that   
 \begin{equation}
  Z_x(s)=V_x(s)-2V_{cx}+\sum_{p\leq cx}\frac{2\log p}{p^{2s}-1}+ O(x^{1/2-\sigma_1})=V_x(s)-2V_{cx}+O(1).
  \end{equation}
Thus we have
$Z_x(s)=V_x(s)-2V_{cx}(s)+O(1)=-W_x(s)+O(1).$\\

Using technique similar to Lemma \ref{lemma2_2} in  part 2, we can show that $Z_x(s)$ has a root at $$s=\sigma_1 +\frac{1}{\log^a x}+O\bigg(\frac{1}{\log^{2a} x}\bigg).$$

\end{proof}
\end{lemma}
\begin{lemma}\label{lemma13_4}
We have the following lower bound on $Z_x(s)$
$$Z_x(s)\gg \frac{x^{1-\sigma_1}}{\log^a x}$$	On $\mathscr{C}_1$.
\begin{proof}
Since the main term of $Z_x(s)$ is just negative of the main term of $W_x(s)$,the proof is identical to the corresponding proof of Lemma \ref{lemma9_2} in part 2.
\end{proof}
\end{lemma}

\section{Proof of Theorem 4}
Here we  state the statement of Theorem 4 for the convenience of the reader.

\begin{theorem}
Let $\sigma_1$ be such that $1/2<\sigma_1<1$ and, let $d_1=\sigma_1-1/2.$
Let $\rho'=\beta'+\di \gamma'$ denote a critical point of the Riemann zeta function such that $\sigma_1<\beta'<1$. Then for any $\epsilon>0$ and for infinitely many $\rho'$ with $\gamma'\to\infty$, we have unconditionally that

\begin{equation}
 \log|\zeta(\rho')|\leq (-B(\s_1)+\epsilon)\frac{(\log\gamma')^{1-\beta'}}{(\log\log\gamma')^{4-2\beta'}},
\end{equation}
where
\be
  B(\sigma_1)  =\frac{(\sigma_1-1/2)^{1-\s_1}\log 2}{(1-\s_1)^2 4^{1-\s_1}}
  \ee
\end{theorem}
\begin{proof}
Let $\mathscr{C}_0$, $\mathscr{C}_1$, $K$ and $\delta$ be the same as defined in part 2.
Let $\tau$ be a real number such that 

$$\big\|{\tau \log p/2\pi}+1/2 \big\|<\delta \quad\hbox{for} \quad p\leq cx,$$ 

 $$\big\|{\tau \log p/2\pi} \big\|<\delta \quad\hbox{for} \quad cx <p\leq x,$$ 
 and $\mathscr{G} $ be the set of such $\tau$ such that above inequalities hold.

We have the following inequality.
\begin{equation}
|V_x(s+\di\tau)-Z_x(s)|\ll \frac{x^{1-\sigma_1}}{\log ^{b-1}x}
\end{equation}
where $s \in \mathscr{C}_1$ and $\tau \in \mathscr{G}.$
The proof is similar to the proof in part 2.

Define $L_x(s)$ as follows.
$$L_x(s)=\sum_{1<n\leq x}\frac{c(n)\Lambda(n)}{n^s \log n}.$$
Clearly,
\begin{equation}\label{eq19_4}
\begin{split}
 \bigg|\sum_{1<n\leq x}\frac{\Lambda(n)n^{-s-\di\tau}}{\log n}-L_x(s)\bigg|
 &= \bigg|\sum_{1<n\leq x}\frac{\Lambda(n)n^{-s-\di\tau}}{\log n}-\sum_{1<n\leq x}\frac{c(n)\Lambda(n)}{n^s \log n}\bigg|\\
&\leq \sum_{1<n\leq x} \frac{\Lambda(n)}{n^\sigma \log n}  | n^{-i\tau} -c(n)|\\
\end{split}
\end{equation}
If $n=p^k$, then $ | n^{-i\tau} -c(n)| = |p^{ik\tau}-c(p)^k|\leq k | p^{i\tau} -c(p)|$. Thus, for $\sigma\geq \sigma_1$
and $\tau \in \mathscr G$, the expression here is
$$
\leq 2\pi \delta \sum_{p^k\leq x}\frac{1}{p^{k\sigma_1}} \ll \delta \frac{x^{1-\sigma_1}}{\log x}
\ll \frac{x^{1-\sigma_1}}{\log^b x}
$$
by our choice of $\delta$ in \eqref{parameters}.
That is, for $\sigma\geq \sigma_1$ and $\tau \in \mathscr G$,
\begin{align}\label{sum L_x}
 \bigg|\sum_{1<n\leq x}\frac{\Lambda(n)}{\log n}n^{-s-\di\tau} -L_x(s)\bigg|
\ll \frac{x^{1-\sigma_1}}{\log^b x}.\end{align}

Also, when $\Re s\geq \sigma_1$,

\begin{align*}
L'_x(s)=-\sum_{n\leq x}\frac{c(n)\Lambda(n)}{n^s}
\ll \sum_{n\leq x}\frac{\Lambda(n)}{n^{\sigma_1}}
\ll x^{1-\sigma_1}.\\
\end{align*}

Thus, if $s$ is on or inside  $\mathscr{C}_1$, then
\begin{equation}\label{eq20_4}
L_x(s)=L_x(\sigma_1)+O\bigg(\frac{x^{1-\sigma_1}}{\log^a x}\bigg). 
\end{equation}

Let's evaluate the value of  $L_x(s)$ at $\sigma_1$.
\label{eq21_4}
\begin{align}
L_x(\sigma_1)&=-2\sum_{p\leq cx}p^{-\sigma_1}+\sum_{p\leq x}p^{-\sigma_1}+O(1)\\
&=\frac{x^{1-\sigma_1}\log c}{(1-\sigma_1)\log^2 x}+O\bigg(\frac{x^{1-\sigma_1}}{\log^{\min(3,a+1)} x}\bigg)
\end{align}
which we get by technique similar in part 2.\\

Thus combining ($\ref{eq19_4}$), ($\ref{eq20_4}$) and ($\ref{eq21_4}$), we see that for $s$ on or inside  $\mathscr C_1$, and $\tau\in \mathscr G$,
 \begin{align}\label{sum L_x}\notag
\Re \sum_{1<n\leq x}\frac{\Lambda(n)}{\log n}n^{-s-\di\tau} =&\frac{x^{1-\sigma_1}\log c}{(1-\sigma_1)\log^2 x}
+O\Big( \frac{x^{1-\sigma_1}}{ \log^{\min(a, 3)} x } \Big)
+O\Big(\frac{x^{1-\sigma_1}}{\log^b x}\Big)\\
=&\frac{x^{1-\sigma_1}\log c}{(1-\sigma_1)\log^2 x}
+O\Big( \frac{x^{1-\sigma_1}}{ \log^{\min(a, 3)} x } \Big).
\end{align}
since $b>a+1$ and $a>2$.

We have the following inequality.
$$\Re\sum_{p^k\leq x}\frac{p^{-k(s+\di\tau)}}{k}\leq \bigg(\frac{\log c}{1-\sigma_1}+o(1)\bigg)\frac{x^{1-\sigma_1}}{\log^2 x}$$
where $s \in \mathscr{C}_1$ and $\tau \in \mathscr{G}.$ 

Let
$$
T_x(s) = \sum_{x<n\leq T^{d_1/4}} \frac{w(n)\Lambda(n)}{\log n} n^{-s}.
$$

Then by Lemma~\ref{lem 12}, when $s\varepsilon\mathscr{C}_1$ and $\tau\varepsilon \mathscr{G}$,
for $s=\sigma+it$, with $\sigma_1 \leq \sigma$ and $|t|\leq 1$, and for $\tau \in R$,  we have 
\begin{equation}\label{ }
\log {\zeta}(s+\di\tau)
=\sum_{1<n\leq x} \frac{\Lambda(n)}{\log n}n^{-s-\di\tau}+T_x(s+i\tau) +O(T^{-d_1^2/16}),    
\end{equation}
and
\begin{equation}\label{eq22_4}
-\frac{\zeta'}{\zeta}(s+\di\tau)
=\sum_{n\leq x}\Lambda(n)n^{-s-\di\tau} -T_{x}^{'}(s+i\tau) +O(T^{-d_1^2/16}).    
\end{equation}

Which we get from \eqref{eq15_2} in part 2,
where $x={d_1\log T}/{4{(\log\log T)}^a}.$ Now we have to show
\begin{equation}\label{eq23_4}
T'_x(s+\di\tau)=  \sum_{x<n\leq T^{d_1/4}}w(n)\Lambda(n)n^{-s-\di\tau}\ll 1,  
\end{equation}
and
\begin{equation}\label{eq24_4}
 T_x(s+\di\tau)=\sum_{x<n\leq T^{d_1/4}}w(n)\Lambda(n)\frac{n^{-s-\di\tau}}{\log n}\ll \frac{1}{\log x}   
\end{equation}
when $ s \in \mathscr{C}_1$ and  for some $\tau \in \mathscr{G}$.\\
By substituting \eqref{eq23_4} in \eqref{eq22_4}, we get
$$-\frac{\zeta'}{\zeta}(s+\di\tau)=Z_x(s)+O\Big(\frac{x^{1-\sigma_1}}{\log^{b-1} x}\Big).$$
Since $Z_x(s)=-W_x(s)+O(1)$, by the same argument of \eqref{W_x 4} and \eqref{zeta'/zeta} we can say that $\zeta'(s)$ has a root for $\rho'$ in $\mathscr{C}_1.$ 
 
We know that
$$ \log c=\log 2/(\sigma_1-1)(1+o(1)).$$
Combining above estimates, we get

\begin{equation}\label{eq14_2}
\log|\zeta(\rho')| \leq  (1+o(1) ) \Big(-\frac{(\sigma_1-1/2)^{1-\s_1}\log 2}{(1-\s_1)^2 4^{1-\s_1}}\Big)
\frac{( \log T)^{1-\sigma_1} }{(\log \log T)^{ 2+b(1-\sigma_1) } }.
\end{equation}

Since $T\leq \gamma'=\Im \rho' \leq 2T$,  \eqref{eq14_2} also holds with $T$ replaced by $\g'$, which is
\eqref{eq1.9}.

Thus we get the following equation.
$$\log|\zeta(\rho')|\leq \big((-B(\s_1)+\epsilon)\big)\frac{(\log\gamma')^{1-\beta'}}{(\log\log\gamma')^{2+b(1-\sigma_1)}},$$ 
where
$$
  B(\sigma_1)  =\frac{(\sigma_1-1/2)^{1-\s_1}\log 2}{(1-\s_1)^2 4^{1-\s_1}}.
  $$

To complete the proof of Theorem~\ref{thm 2} it  only  remains to show the existence of a $\tau \in R\cap\mathscr G$ satisfying \eqref{eq23_4} and \eqref{eq24_4}.
We can apply the same idea and Kronecker's formula with minor changes. The proof will be exactly the same as the proof in part 2 with  a different set of $\beta_p$.  

Combining all the estimates and choosing $b=3+\epsilon''$, we get
$$\log|\zeta(\rho')|\leq \big((-B(\s_1)+\epsilon)\big)\frac{(\log\gamma')^{1-\beta'}}{(\log\log\gamma')^{5-3\beta+\epsilon'}},$$ 
where $\epsilon'=\epsilon''(1-\sigma).$

Thus our proof is complete.
\end{proof}

\section{Acknowledgments}
{
The author gives sincere thanks to his doctoral advisor Steven M. Gonek for
introducing the problem in this paper and also for providing guidance and support
during the process of its study. Professor Gonek also read an earlier version of
this paper and made many useful suggestions which significantly improved the
exposition.

}


\begin{thebibliography}{9}

\bibitem{paper}
Steve M. Gonek and Hugh L. Montgomery,Extreme values of zeta function at critical points, \textit{Q. J. Math.} \textbf{67} (2016), no. 3, 483–505.


\bibitem{paper3}
J.E. Littlewood, On the Riemann zeta function, \textit{Proc. London math. Soc.}(2) \textbf{24} (1926), 175-201.

\bibitem{paper4}
J.E. Littlewood, On the function $1/\zeta(1+\di t)$, \textit{Proc. London math. Soc.}(2) \textbf{27} (1928), 349-357.

\bibitem{Mont-Thomp}
Montgomery H L and Thompson J G,   Geometric properties of the zeta function \textit{Acta Arith.} 
\textbf{155} (2012), 373–96.

\bibitem{paper2}
Hugh L. Montgomery and R.C. Vaughan,Hilbert's inequality, \textit{J. London Math Soc.}(2) \textbf{8} (1974), 73-82.

\bibitem{classical} 
	Hugh L. Montgomery and Robert C. Vaughan,
	\textit{Multiplicative Number Theory 1: Classical Theory}, 
	Cambridge University Press, Cambridge,2007.



\bibitem{paper5}
E.C.Tichmarsh, On the function $1/\zeta(1+\di t)$,\textit{Quart. J. math.(Oxford)},\textbf{4} (1933),64-70.

\bibitem{Titch}  E.C. Titchmarsh, The theory of the Riemann
zeta-function (2nd edition),  Oxford University Press, Oxford, 1986.

\end{thebibliography}
\end{document}